# The Generalized Burnside Theorem
## S. Bachmuth

## 1. Introduction

We use the same notation that was used in [1]. Unless specifically mentioned otherwise, all groups are 2-generator. We fix a prime power $q = p^e$. There is no restriction on q as long as it is a prime power.

The four matrix groups $F(\mathcal{R}[t, t^{-1}])$, $F(\mathcal{R})$, $F(S[t, t^{-1}])$ and $F(S)$ were introduced in [1] and are defined in Section 2. A proof of any claim about these groups not found in this paper can be found in [1]. We list here some relevant facts concerning these groups:

$\mathcal{R} = Z[x, x^{-1}, y, y^{-1}]$ is the 2 variable Laurent polynomial ring with integer coefficients.
$\mathcal{R}[t, t^{-1}] = Z[x, x^{-1}, y, y^{-1}, t, t^{-1}]$ is the 3 variable ring.
$F(\mathcal{R}[t, t^{-1}])$ is a free group of 2x2 matrices with entries in $\mathcal{R}[t, t^{-1}]$.
$F(\mathcal{R})$ is a free metabelian group of 2x2 matrices with entries in $\mathcal{R}$. It is the image of $F(\mathcal{R}[t, t^{-1}])$ upon setting $t = 1$ in each matrix.
For each prime power q, there exists a quotient ring $S = S(q)$ of $\mathcal{R}$ such that upon replacing the elements of $\mathcal{R}$ by those of $S$ in $F(\mathcal{R})$ and also in $F(\mathcal{R}[t, t^{-1}])$, we have
1) $F(S)$ is a free metabelian exponent q group of 2x2 matrices over $S$. That is, $F(S)$ is isomorphic to $F/F''F^q$, where F is a free group.
2) $F(S[t, t^{-1}])$ is a solvable group of 2x2 matrices over $S[t, t^{-1}]$.

We continue to reserve $\underline{\alpha}$ for the map
$$\underline{\alpha}: F(\mathcal{R}) \to F(S)$$
which arises from the ring homomorphism $\mathcal{R} \to S$. Thus, $\underline{\alpha}$ is the Burnside metabelian map. (The Burnside map is $\alpha: F(\mathcal{R}[t, t^{-1}]) \to G$, where $G \cong F/F^q$.)

Our goal is to establish

<u>Theorem</u> GB (Generalized Burnside): Let γ be a mapping of the free group $F(\mathcal{R}[t, t^{-1}])$ onto a group G such that γ induces $\underline{\alpha}$. Then G (i.e., the image of γ) is solvable.

We will review the definition and relevant facts about induced maps in Section 3. The assumption that γ induces $\underline{\alpha}$ implies many things. What is not at all obvious is that it implies that, in addition to having at least one generator of order q, the commutator subgroup of G has exponent q. The presence of this large subgroup of exponent q must surely be the reason for the identities which imply solvability.

Since the Burnside map on the free group $F(\mathcal{R}[t, t^{-1}])$ induces $\underline{\alpha}$, we have the following:

Corollary: The Burnside group of exponent q is solvable.

Among the solvable groups arising from Theorem GB, some are finite (e.g., the Burnside groups of prime-power exponent), while others possess elements of infinite order. $F(S[t, t^{-1}])$ is an example of the latter. But whether finite or infinite, all of these groups are closely associated with exponent q groups. Indeed, each has a subgroup of exponent q which is strictly larger than the commutator subgroup.

Theorem B in [2] allows good estimates of upper bounds for the solvability class of Burnside groups for 2 generators (and more generally for $k \leq p+1$ generators). For example, a straightforward calculation will show that the 2 generator Burnside group of exponent $q = p^e$ has (solvability) class at most n where $2^n \geq e(p^e - p^{e-1}) + 1$. Moreover, good lower bounds can easily be obtained by substituting primitive roots of unity for x and y in the ring $S$. It appears likely that these upper and lower bounds are not far apart, and in some cases the same for the 2 generator groups.

The proof of the Generalized Burnside Theorem is broken up into two parts.

Proposition 1: Suppose $\gamma : F(\mathcal{R}[t, t^{-1}]) \to G$ induces $\underline{\alpha}$. Then the entries of $F(\mathcal{R}[t, t^{-1}])$ which are in $\mathcal{R}$ are sent into S.

Proposition 2: Suppose the surjective map $\gamma : F(\mathcal{R}[t, t^{-1}]) \to G$ sends $\mathcal{R}$ into $S$. Then G (the image of γ) is solvable.

As we propose to make clear, both of these propositions were in essence proved in [1]. The proof of Proposition 1 given here is similar to the arguments given in Section 5 of [1]. The format for Proposition 2, which solves a problem in groups (e.g., establishing solvability) via calculations in a ring, goes back at least to W. Magnus' classic paper [3]. In that paper and numerous others by various authors since, the terms in a higher commutator series are linked to a 'degree function' in a ring - in our case to the powers of an augmentation ideal. Proposition 2 is achieved by showing that as one proceeds down the derived series of $F(\mathcal{R}[t, t^{-1}])$, the powers of the augmentation ideal of $\mathcal{R}$ increase. This means that Theorem B in the joint paper with Heilbronn and Mochizuki [2] ensures that they eventually fall into the q-cyclotomic ideal of $\mathcal{R}$ (i.e., the zero of $S$). Thus, upon the transfer of $F(\mathcal{R}[t, t^{-1}]) \to G$, the assumption $\mathcal{R} \to S$ implies that G is solvable.

The following sections will review the necessary background and then present the proof of Propositions 1 and 2, and thus Theorem GB. .

We are grateful to Bob Guralnick for pointing out the error in the original abstract. I take pleasure in thanking him for his helpful suggestions.

## 2. Required Background

In reviewing material from [1] in this section, we give the formulation for k generators. The proof of Theorem GB, however, will be limited to 2 generators. Additional details of the concepts may be found in [1].

$\mathcal{R} = \mathcal{R}(k) = Z[x_1, x_1^{-1}, \ldots, x_k, x_k^{-1}]$ is the k generator Laurent polynomial ring over the integers. $\mathcal{R}[t, t^{-1}] = Z[x_1, x_1^{-1}, \ldots, x_k, x_k^{-1}, t, t^{-1}]$ is the k+1 generator ring. $\mathcal{I}(q)$ is the q-cyclotomic ideal in $\mathcal{R}$, the ideal generated by all q-cyclotomic elements $1 + u + u^2 + \cdots + u^{q-1}$ where u is a positive unit in $\mathcal{R}$ and q is the prime power $q = p^e$. $S$ is the quotient ring $\mathcal{R}/\mathcal{I}(q)\Sigma$, where $\Sigma$ is the augmentation ideal of $\mathcal{R}$.

We can now describe the k x k matrix groups $F(\mathcal{R})$, $F(S)$, $F(\mathcal{R}[t, t^{-1}])$ and $F(S[t, t^{-1}])$.

Let v be the k-tuple, $v = (1-x_1, \ldots, 1-x_k)$ and N the k x k matrix $N = [\lambda_i v]$, where $\lambda_i v$ denotes the ith row of N, $i = 1,\ldots,k$. $F(\mathcal{R})$ is the set of all matrices $uI + N$, where I is the identity matrix, u is a (positive) unit in $\mathcal{R}$, and the $\lambda_1, \ldots, \lambda_k$ in $\mathcal{R}$ satisfy $\lambda_1(1-x_1) + \ldots + \lambda_k(1-x_k) = 1-u$. For a matrix M in $F(\mathcal{R})$ we often omit the I and write $M = u + N$. Alternatively, one may define $F(\mathcal{R})$ by the generators:
$$F(\mathcal{R}) = gp\langle M_1, \ldots, M_k \rangle, \text{ where for } j = 1, \ldots, k$$
$M_j = x_j I + [\lambda_i v]$, $\lambda_i = 0$ unless $i = j$ in which case $\lambda_j = 1$.

Replacing $\mathcal{R}$ by $S$ defines $F(S)$ as well as $F(S[t, t^{-1}])$ as soon as we define $F(\mathcal{R}[t, t^{-1}])$.

In defining $F(\mathcal{R}[t, t^{-1}])$ we first define k x k matrices $T_i$ ($i = 2,\ldots,k$). $T_i$ has t in the first i-1 diagonal entries, 1 in the remaining diagonal entries, 1-t in the ith row prior to the diagonal, and zeros everywhere else. Notice that $T_2 = T$ in [1]. It is easy to see that $gp\langle T_2, \ldots, T_k \rangle$ is free abelian of rank k-1, but all that is necessary is that it be abelian (in generalizing Lemma 5 of [1] for arbitrary ranks). We can now define $F(\mathcal{R}[t, t^{-1}])$.
$$F(\mathcal{R}[t, t^{-1}]) = gp\langle M_1, M_2 T_2, M_3 T_3, \ldots, M_k T_k \rangle$$

Those wishing to deal with the Burnside problem for k generators can now follow the procedures of [1] or those that will follow. The proof that $F(\mathcal{R})$ is free metabelian of rank k and that $F(S)$ is the free metabelian Burnside group of exponent q and rank k should follow as in [1]. (That $F(\mathcal{R})$ and $F(S)$ are images of $F(\mathcal{R}[t, t^{-1}])$ and $F(S[t, t^{-1}])$ respectively is obvious.) As in [1], we have the commutative square where the horizontal maps $\alpha$ and $\underline{\alpha}$ come from the ring homomorphism $\mathcal{R} \to S$ and the vertical maps are those determined by sending t to the identity.

$$\begin{array}{ccc} F(\mathcal{R}[t, t^{-1}]) & \xrightarrow{\alpha} & F(S[t, t^{-1}]) \\ \downarrow & & \downarrow \\ F(\mathcal{R}) & \xrightarrow{\underline{\alpha}} & F(S) \end{array}$$

These preparations should be helpful for those wishing to delve into the k-generator groups, where much work remains in the determination of bounds. From this point forward we will deal only with 2-generator groups. When $k = 2$, our notation reduces to that of [1]. Hence for the remainder of the paper,

$\mathcal{R} = \mathcal{R}(2) = Z[x, x^{-1}, y, y^{-1}]$ and $F(\mathcal{R})$, $F(S)$, and $F(\mathcal{R}[t, t^{-1}])$ are as in [1].
$F(\mathcal{R}) = gp\langle M_1, M_2 \rangle$ and $F(\mathcal{R}[t, t^{-1}]) = gp\langle M_1, M_2T \rangle$.

## 3. Proof of the Generalized Burnside Theorem

<u>Theorem</u> GB (Generalized Burnside): Let $\gamma$ be a mapping of the free group $F(\mathcal{R}[t, t^{-1}])$ onto a group G such that $\gamma$ induces $\underline{\alpha}$. Then G (i.e., the image of $\gamma$) is solvable.

The proof is a consequence of the following two propositions:

<u>Proposition 1</u>  Suppose $\gamma : F(\mathcal{R}[t, t^{-1}]) \to G$ induces $\underline{\alpha}$. Then the entries of $F(\mathcal{R}[t, t^{-1}])$ which are in $\mathcal{R}$ are sent into S.

<u>Proposition 2</u>  Suppose the surjective map $\gamma : F(\mathcal{R}[t, t^{-1}]) \to G$ sends $\mathcal{R}$ into S. Then G (the image of $\gamma$) is solvable.

Before proving Proposition 1, we first review the concept of induced map. Refer to Section 5 of [1] for a more detailed discussion.

For g in $F(\mathcal{R}[t, t^{-1}])$, let g/t be the matrix obtained by setting $t = 1$ in g. Thus g/t is in $F(\mathcal{R})$.

<u>Def</u>: Given the map $\lambda: F(\mathcal{R}[t, t^{-1}]) \to G$. Define $\underline{\lambda}: F(\mathcal{R}) \to G/G''$ by $\underline{\lambda}(g/t) = \lambda(g)G''$ for g in $F(\mathcal{R}[t, t^{-1}])$. $\underline{\lambda}$ is the map induced from $\lambda$.

That $\underline{\lambda}$ is well defined is just Lemma 9 of [1] found in Section 5 of [1]. The proof depends on the fact that $\lambda$ takes the second derived group of $F(\mathcal{R}[t, t^{-1}])$ into $G''$.

<u>Remark</u>: More generally, as is well known, one can define an induced map on any quotient of a free group by a fully invariant subgroup. For example, questions concerning automorphisms induced from free groups have been studied by various authors in a variety of allowable quotient groups. Our concern will focus on maps of

$F(\mathcal{R}[t, t^{-1}])$ which induce $\underline{\alpha}$, the Burnside metabelian map.

Proposition 1.  Suppose $\gamma : F(\mathcal{R}[t, t^{-1}]) \to G$ induces $\underline{\alpha}$. Then $\gamma$ sends the entries of $\mathcal{R}$ into S.

Proof: This is actually Lemma 11 in [1], but we shall give a proof again. The proof begins by noting that since $M_1$ does not have t in any entry, $M_1/t = M_1$. Thus, since $M_1$ is in both $F(\mathcal{R})$ and $F(\mathcal{R}[t, t^{-1}]$, we have $\gamma(M_1/t) = \gamma(M_1) = \underline{\gamma}(M_1)$. By hypothesis, $\underline{\gamma}(M_1) = \underline{\alpha}(M_1)$, and by definition, $\underline{\alpha}(M_1)$ takes the entries of $M_1$ into S. However, the entries of $M_1$ contain the generators of $\mathcal{R}$ and thus $\gamma$ sends all entries from $\mathcal{R}$ to S.

Comments: This Proposition already tells us much about G.
1)  The image of the first generator of $F(\mathcal{R}[t, t^{-1}])$ is the same as the image of the first generator of the metabelian Burnside map and is therefore an element of order q.
2)  Since the elements of of $\mathcal{R}$ are sent into S., G must be a matrix group over a quotient of a Laurent polynomial ring.
3)  Thus far, $\gamma$ is behaving exactly like the map of $F(\mathcal{R}[t, t^{-1}])$ to $F(S[t, t^{-1}])$ studied in [1]. All that remains to define $\gamma$ is to know what happens to t. But, at this point we can be sure that whatever happens to t, the ideal defining G will contain the ideal defining $F(S[t, t^{-1}])$. Thus G will be a homomorphic image of the (solvable) group $F(S[t, t^{-1}])$. This is actually the first proof of Proposition 2 coming up.

Note: Theorem B in [2], which is used in Proposition 2, links the augmentation and cyclotomic ideals of $\mathcal{R}$. It is summarized in Section 2 of [1] (Lemma 3 (i) ).

Proposition 2: Suppose the surjective map $\gamma : F(\mathcal{R}[t, t^{-1}]) \to G$ sends $\mathcal{R}$ into S. Then G is solvable.

Proof: First note that for $G = F(S[t, t^{-1}])$, Proposition 2 has already been established. That is, for this choice of G, we have the surjective map $\gamma: F(\mathcal{R}[t, t^{-1}]) \to F(S[t, t^{-1}])$ in which $\gamma$ sends $\mathcal{R}$ to $S$ while leaving t unchanged. We may therefore invoke Lemma 8 in [1] to conclude that $F(S[t, t^{-1}])$ is solvable. We complete the proof by noting that G in Proposition 2 is an (homomorphic) image of $F(S[t, t^{-1}])$. Namely, by assumption, the map $F(\mathcal{R}[t, t^{-1}]) \to G$ sends $\mathcal{R}$ to $S$. This means that G is constructed by means of a quotient of an ideal of $\mathcal{R}[t, t^{-1}]$ that contains the ideal $\mathcal{I}(q)\Sigma$ used to define S. But, the ideal of $\mathcal{R}[t, t^{-1}]$ used in the construction of $F(S[t, t^{-1}])$ is exactly $\mathcal{I}(q)\Sigma$ and thus smaller or possibly equal to the ideal used to construct G. Thus G is an image of $F(S[t, t^{-1}])$.

Alternate Proof: We again begin by observing that for $G = F(S[t, t^{-1}])$, Proposition 2

has already been established by Lemmas 4 through 8 in [1]. For this choice of G, the map $F(\mathcal{R}[t, t^{-1}]) \to F(S[t, t^{-1}])$ takes $\mathcal{R}$ to $S$ and leaves t fixed. Consider now arbitrary G. By hypothesis, the elements of $\mathcal{R}$ are sent to $S$ in the map $\gamma: F(\mathcal{R}[t, t^{-1}]) \to G$. Although we no longer have information about t, we can still invoke Lemmas 4 through 8 in [1] because only the augmentation ideal of $\mathcal{R}$ is involved. Indeed, in the proofs of Lemmas 4 through 7, the computation of the derived series is done in the free group $F(\mathcal{R}[t, t^{-1}])$. Only in Lemma 8 does $F(S[t, t^{-1}]$ appear where the information from $F(\mathcal{R}[t, t^{-1}]$ is transferred to $F(S[t, t^{-1}])$ via the correspondence $\mathcal{R} \to S$. G replaces $F(S[t, t^{-1}]$, but the proof of solvability remains the same. Since $\gamma$ sends the elements of $\mathcal{R}$ to $S$ in the transfer from $F(\mathcal{R}[t, t^{-1}])$ to G, Theorem B in [2] assures that the elements in a large enough power of the augmentation ideal of $\mathcal{R}$ are sent to the zero of S. This means that in the expansion of the derived series in G, the elements in a high enough solvability class become the identity in G.

Remarks:
1) In this second proof of Proposition 2, we are in effect leaving intact the procedure of Vaughan-Lee as detailed in [1] which showed that $F(S[t, t^{-1}])$ is solvable. We use the hypothesis $\mathcal{R} \to S$ to show that the image of the map $F(\mathcal{R}[t, t^{-1}]) \to G$ is solvable in place of the original map $F(\mathcal{R}[t, t^{-1}]) \to F(S[t, t^{-1}])$, a map which was defined by $\mathcal{R} \to S$.

2) The image of t determines, among other things, whether G is finite or infinite and is also important in the determination of the actual solvability class. If t is left fixed, then $G = F(S[t, t^{-1}])$ is not only an infinite solvable group, but also has the largest possible solvability class for the chosen q. Enlarging the ideal will generally decrease the solvability class. In the other extreme case of sending t to 1, the ideal becomes as large as possible and G, which becomes $F(S)$, has the smallest possible solvability class, (that is, metabelian or abelian when q = 2).

3) The commutator subgroup of $F(S[t, t^{-1}])$ has exponent q (see Section 4 of [1]. Thus all it's non abelian images are also exponent q. Because these GB groups have at least one generator of order q, their commutator subgroup is a proper subgroup of their exponent q subgroup Thus these generalized Burnside groups do not stray too far from exponent q groups. Surely these huge numbers of exponent q elements yield the relations which imply solvability. Unfortunately, these and related identities are known only for very small q. Even q = 5 has thus far remained a complete mystery.

Caution: Those interested in k-generators should realize that Theorem B in [2] is valid only for $k \leq p+1$. In [2] it is shown that the bound $e(p^e - p^{e-1})$ increases when $k > p+1$. For a proof of solvability it is not difficult to get around this, but at the price of losing control of the solvability class. An extension of Theorem B for $k > p+1$ would greatly enhance our knowledge of bounds.

# References

[1]   S. Bachmuth,  Solvable matrix groups and the Burnside Problem, ArXiv: math.GR/0609415,   14 Sept. 2006, 14 pages

[2]   S. Bachmuth, H. Heilbronn and H. Y. Mochizuki,  Burnside Metabelian groups, Proc. Royal Soc. (London) 307 (1968), 235-250.
      Reprinted in The Collected Papers of Hans Arnold Heilbronn, E.J. Kani and R.A. Smith editors, J. Wiley and Sons.

[ 3]   W. Magnus, Beziehungen zwischen Gruppen und Idealen in einem speziellen Ring. Math. Ann., 111, 259-280 (1935).